%% file: S-Galerkin-a.tex
\documentclass[a4paper, notitlepage, 12pt]{article}

\usepackage{ucs}                
\usepackage[utf8x]{inputenc} 	

\usepackage[T1]{fontenc}     
\usepackage{lmodern}         

\usepackage{a4wide} 		    
\usepackage{color}
\usepackage{lscape}
\usepackage{amsmath}
\usepackage{amssymb}
\usepackage{amsbsy}
\usepackage{amsthm}
\usepackage{wasysym}
\usepackage{float}
\usepackage{algpseudocode}
\usepackage[section]{algorithm}
\usepackage[breaklinks]{hyperref}
\usepackage{authblk}

\usepackage{multirow} 
\usepackage{array}
\usepackage{graphicx}

\usepackage[british]{babel}

\usepackage{refdef}
\usepackage{nfontdef}

\include{nformula}

\newcommand{\ip}[2]{\langle #1, #2 \rangle}

\definecolor{myred}{rgb}{1, 0.2, 0.2}

\newcommand{\authorhgm}{\authorcr Hermann G. Matthies}
\newcommand{\authordl}{Dishi Liu}
\newcommand{\authoral}{Alexander Litvinenko}
\newcommand{\authorlg}{Lo\"{i}c Giraldi}
\newcommand{\authoran}{Anthony Nouy}

\newcommand{\affilecn}{\'Ecole Centrale de Nantes, GeM, UMR CNRS 6183, LUNAM Universit\'e}
\newcommand{\affilwire}{Institute of Scientific Computing, 
                        Technische Universit\"at Braunschweig}
\newcommand{\affildlr}{C$^2$A$^2$S$^2$E, German Aerospace Center (DLR),
                        Braunschweig}

\newcommand{\thetitle}{To be or not to be intrusive? \authorcr
       The solution of parametric and stochastic equations\authorcr
       --- the ``plain vanilla'' Galerkin case}
\newcommand{\theauthor}{\authorlg, \authoral, \authordl, \authorhgm, \authoran}
\newcommand{\thesubject}{(MSC) \tbf{65M70}, 65J15, 60H25, 65D30, 34B08, 35B30}
\newcommand{\thekeywords}{parametric problems, stochastic equation, uncertainty
            quantification, Galerkin approximation, coupled system, 
            non-intrusive computation}
\newcommand{\textdate}{\today}

\newcommand{\thebib}{./bib}



\begin{document}

\title{\thetitle\thanks{Work partly supported by the Deutsche
          Forschungsgemeinschaft (DFG) through the SFB 880.}}
\author[c]{\authorlg}
\author[a]{\authoral}
\author[b]{\authordl}
\makeatletter
\author[a]{\authorhgm 
\thanks{Corresponding author: TU Braunschweig, D-38092 Braunschweig, 
       Germany, e-mail: \texttt{wire@tu-bs.de}}
}
\makeatother
\author[c]{\authoran}

\affil[a]{\affilwire}
\affil[b]{\affildlr}
\affil[c]{\affilecn}

\date{\today}


\maketitle

\input{\thetext/abstract}








\input{\thetext/introduction}
\input{\thetext/parametric-gal}
\input{\thetext/galerkin}
\input{\thetext/residual}
\input{\thetext/iteration}
\input{\thetext/conclusion}



\bibliography{\thebib/jabbrevlong,\thebib/matthies_BU_paper-1,\thebib/phys_D,\thebib/fa,\thebib/risk,\thebib/fuq-new,\thebib/highdim}


{ 
   \tiny
   \begin{verbatim}
    $Id: S-Galerkin-a.tex,v 1.10 2013/09/05 09:51:00 hgm Exp $
   \end{verbatim}
   }



\end{document}

%% file: nformula.tex
\newtheorem{thm}{Theorem}

\newtheorem{prop}[thm]{Proposition}

\newtheorem{coro}[thm]{Corollary}

\newtheorem{lem}[thm]{Lemma}

\newcommand{\ignore}[1]{}

\newcommand{\vrho}{\varrho}
\newcommand{\vepsilon}{\varepsilon}

\newcommand{\vphi}{\varphi}

\newcommand{\vek}[1]{\mathchoice{\displaystyle\boldsymbol{#1}}
{\textstyle\boldsymbol{#1}}{\scriptstyle\boldsymbol{#1}}
{\scriptscriptstyle\boldsymbol{#1}}}
\newcommand{\mat}[1]{\mathchoice{\displaystyle\mathbf{#1}}
{\textstyle\mathbf{#1}}{\scriptstyle\mathbf{#1}}
{\scriptscriptstyle\mathbf{#1}}}

\newcommand{\Mtil}[1]{\mat{\tilde{#1}}}

\newcommand{\EXP}[1]{\mathbb{E}\left(#1\right)}

\newcommand{\spn}{\mathop{\mathrm{span}}\nolimits}

\newcommand{\di}{\mathrm{d}}
\newcommand{\Di}{\mathrm{D}}

%% file: chapters/abstract.tex
%

\begin{abstract}
In parametric equations---stochastic equations are a special case---one may
want to approximate the solution such that it is easy to evaluate its
dependence of the parameters.  Interpolation in the parameters is an obvious
possibility, in this context often labeled as a collocation method.
In the frequent situation where one has a ``solver''
for the equation for a given parameter value---this may be a software component
or a program---it is evident that this can independently solve for
the parameter values to be interpolated.
Such \emph{uncoupled} methods which allow the use of the original solver are
classed as ``non-intrusive''.  By extension, all other
methods which produce some kind of \emph{coupled} system are often---in
our view prematurely---classed as ``intrusive''.  
We show for simple Galerkin formulations
of the parametric problem---which generally produce coupled systems---how one
may compute the approximation in a \emph{non-intusive} way.

\vspace{5mm}
{\noindent\textbf{Keywords:} \thekeywords}

\vspace{5mm}
{\noindent\textbf{Classification:} \thesubject}

\end{abstract}

%
%
%
%
%
%
%
%


%% file: chapters/introduction.tex
%

\section{Introduction}  \label{S:intro}
Many problems depend on parameters, which may be a finite set of
numerical values, or mathematically more complicated objects like for
example processes or fields.  We address the situation where we have
an equation which depends on parameters; stochastic equations are a
special case of such parametric problems where the parameters are
elements from a probability space.  One common way to represent this
dependability on parameters is by evaluating the state (or solution) of
the system under investigation for different values of the parameters.
Particularly in the stochastic context this ``sampling'' is a common
procedure.  But often one wants to evaluate the solution quickly for a
new set of parameters where is has not been sampled.  In this
situation it may be advantageous to express the parameter dependent
solution with an approximation which allows for rapid evaluation of
the solution or functionals thereof---so called quantities of interest
(QoI)---in dependence of the parameters.  Such approximations are also
called \emph{proxy} or \emph{surrogate} models, \emph{response functions},
or \emph{emulators}.  This last term was chosen so as to contrast with
\emph{simulator}, which is the original solver for the full equation.
Such approximations are used in several fields, notably optimisation
and uncertainty quantification, where in the last case the parameters
are random variables and one deals with stochastic equations.  All
these methods may be seen as \emph{functional approximations} --- 
representations of the solution by an ``easily computable'' function
of the parameters, as opposed to pure \emph{samples}.

The most obvious methods of approximation used are based on
interpolation, in this context often labelled as \emph{collocation
  methods}.  In this case it is usually assumed that the parameters
are in some sub-domain of a manifold, usually simply just in some
finite-dimensional vector space, and the interpolation is often on
\emph{sparse grids} \cite{barthelmannEtAl99, babuska2007stochastic,
  gana2005zaba, xiu2005highorder}.  This process normally gives the
approximation (interpolant) as a finite linear combination of some
basis functions used for the interpolation, often global multi-variate
polynomials \cite{xiuKarniadakis02a}, or piecewise polynomials
\cite{babuska2004galerkin, xwan2005MEgPC}, or methods based on radial
basis functions, kriging, or neural networks.

Another approach is to again choose a similar finite set of basis
functions, but rather than interpolation use some other projection 
onto the subspace spanned by these functions.  Usually this will
involve minimising some norm of the difference between the true
parametric solution and the approximation, and in many cases this norm
will be induced by an inner product, often in the form of an integral
w.r.t.\ some measure---in the case of stochastic equations this will
be the underlying probability measure.  These integrals in turn may be
numerically evaluated through quadrature formulas---often again on
sparse \emph{Smolyak} or adaptive grids \cite{novakRitter96-highdim,
  gerstnerGriebel98-numint, petras01-fast, keeseMatthies03siam,
  gerstner2003dimensionadaptive}---which need evaluations of the
integrand---part of which is the parametric solution---at a finite
number of parameter values.  Such methods are sometimes called
\emph{pseudo-spectral projections}, or \emph{regression solutions}, or
\emph{discrete projections} \cite{cohen2013onthestability,
  constantine2012sparse, migliorati2012analysis, clenet2012,
  poette2012non, blatman2011adaptive, hosder2010pointcollocation,
  berveiller2006stochastic, reagan2003uncertainty}.

In the frequent situation where one has a ``solver'' for the equation
for a given parameter value, i.e.\ a software component or a program,
it is evident that this can be used to independently---i.e.\ if
desired in parallel---solve for all the parameter values which
subsequently may be used either for the interpolation or in the
quadrature for the projection.  Such methods are therefore
\emph{uncoupled} for each parameter value, and obviously allow to use
the original solver.  Therefore, they additionally often carry the
label ``non-intrusive''.  Without much argument all other
methods---which produce a coupled system of equations---are almost
always labelled as ``intrusive'', meaning that one cannot use the
original solver, e.g.\ \cite{eldred2008evaluation, herzog2008intrusive,
stefanou2009stochastic, eldred2011design, DoostanAavGi2013}.  We
want to show here that this not necessarily the case.

Like most methods which are based on the solution at discrete
parameter values, the non-intrusive methods mentioned above ``forget''
the original equation, i.e.\ the fact that the approximation has to
satisfy the parametric equation.  This is generally the state of
affairs when using the proxy model in the domain of optimisation.  On
the other hand, methods which try to ensure that the approximation
satisfies the parametric equation as well as possible are often based
on a Rayleigh-Ritz or Galerkin type of ``ansatz'', which leads to a
\emph{coupled} system for the unknown coefficients
\cite{ghanemSpanos91, matthies:1999:CMAME, xiuEtAl02c,
  babuska2004galerkin, matthiesKeese05cmame, frSchwab05, TodorSCA,
  Matthies-2008-ZAMM, cohen2011analytic}.  This is often taken as an
indication that the original solver can not be used, i.e.\ that these
methods are ``intrusive''.  But in many circumstances these methods
may as well be used in a \emph{non-intrusive} fashion.  Although there
are some publications concerning special cases of non-intrusive
Galerkin-like methods \cite{acharjee2007non, constantine2009hybrid,
  constantine2010spectral}, this has not been widely recognised as a
general possibility.  A kind of in-between possibility is the
so-called \emph{reduced basis method}, see \cite{BoyLeBYMadE:2010,
  buffa2012priori} for recent expositions.  Here a new basis for the
parametric solution is built from solves at particular parameter
values, but the ``interpolation'' is achieved by a Galerkin projection
onto the spanned subspace.  This method also establishes a connection
between proxy models and reduced order models, something we will not
pursue further here.

Recent developments for \emph{low-rank separated} approximations
\cite{Falco:2012, chinesta2011short, Khoromskij:2011, NouyACM:2010,
 Doostan2009, Nouy2009, Espig:2013} of parametric or stochastic equations
are based on the minimisation of a least squares or similar
functional, and naturally lead to Galerkin-type equations.  Although
it is important to show that these can also be dealt with in a
non-intrusive manner, here we concentrate on the ``plain vanilla'',
i.e.\ standard, Galerkin case.  Non-intrusive computation of separated
approximations will be investigated elsewhere.

Most of the literature cited so far is concerned with the case of
stochastic equations, and although these are a special case of
parametric equations, the methods and techniques used there may be
used in the wider context of general parametric equations, see
\cite{boulder:2011} for a synopsis of these connections of such
parametric problems.

The question whether a method is intrusive or not is often very
important in practise.  The ``solver'' (for a single parameter value)
may contain much specialised knowledge, and may therefore represent
quite a valuable investment of effort.  In case the method is labelled
intrusive, it may seem like the whole---often very domain
specific---process and effort of producing a solver, now for the
coupled Galerkin system, would have to be repeated again.  Therefore,
in many cases the wish to re-use existing software guides the choice
of method.  But as already mentioned, some very effective new methods
based on low-rank approximations fall in the class of ``not obviously
non-intrusive'' methods; hence as a first step we show here that the
simple ``plain vanilla'' coupled Galerkin method may be computed
non-intrusively, the low-rank approximation case will be treated
elsewhere.

A method for a parametric problem will be here considered intrusive if
one has to modify the original software to solve the parametric
problem.  Thus it turns out that the question of whether a method is
intrusive or not hinges on what kind of interface one has to the
software, and is thus a software-engineering question.  Most often it
is possible to not only compute the solution for a certain parameter
value---the solver being usually iterative---but also the residuum or
a ``preconditioned residuum'' given a ``trial solution''.  This
usually means---for the preconditioned residuum---doing just one
iteration with the solver instead of iterating all the way to
convergence.  This is the kind of interface which will be assumed
here, and we shall show that this can be used without any change to
solve the Galerkin equation.

The plan for the rest of the paper is as follows: In the following
\refS{param} we introduce the notation and assumptions for the
parametric problem.  In \refS{galerkin} we introduce the Galerkin
approximation, describe alternative formulations, and prove the
convergence and speed of a basic \emph{block-Jacobi} algorithm for the
coupled Bubnov-Galerkin system.  In the \refS{resid} it is shown how
the residual in the iteration may be computed \emph{non-intrusively},
mainly via numerical integration.  The behaviour of the modified
iterates is analysed, and it is shown that they accumulate in the
vicinity of the solution.  A small numerical example is investigated
in \refS{iterate}, it shows how the non-intrusive computation works,
and it confirms the theoretical predictions.

%
%
%
%
%
%
%
%
%
%
%


%% file: chapters/parametric-gal.tex
%

\section{Parametric Problems} \label{S:param}
To be more specific, let us consider the following situation: we are
investigating some physical system which is modelled by an equation
for its state $u \in \C{U}$ --- a Hilbert space for the sake of
simplicity,
\begin{equation} \label{eq:I}
   A(p;u) = f(p),
\end{equation}
where $A$ is an operator modelling the physics of the system, and
$f\in\C{U}^*$ is some external influence (action / excitation /
loading).  The model depends on some parameters $p \in \C{P}$.  In
many cases \refeq{eq:I} is the abstract formulation of a partial
differential equation.  But for the sake of simplicity we shall assume
here that we are dealing with a model on a finite-dimensional space
$\C{U}$ with $N := \dim \C{U}$, e.g.\ a partial differential equation
after discretisation.  For simplicity we will identify $\C{U}$ and
$\C{U}^*$, and if needed we will use an orthonormal basis
$\{v_n\}_{n=1}^N$, i.e.\ $\spn \{v_n\}_{n=1}^N = \C{U}$ and
$\ip{v_n}{v_m}_{\C{U}} = \delta_{n,m}$, the Kronecker-$\delta$.

Assume that for all $p\in\C{P}$, \refeq{eq:I} is a well-posed problem.
This means that $A$ as a mapping $u \mapsto A(p;u)$ for a fixed $p$ is
bijective and continuously invertible, i.e.\ for each $p$
and $f$ it has a unique solution, which will be denoted by $u^*(p)$,
such that for all $p: \; A(p;u^*(p)) = f(p)$.

Although this will not be needed here, let us remark that if the map
$A$ were also differentiable w.r.t.\ $u$, well-posedness implies that
this partial derivative $\Di_u A$ is non-singular and also
continuously invertible.  Now --- if the set $\C{P}$ has a
differentiable structure, e.g.\ if it is a differentiable manifold or
even a vector space --- one may invoke a version of
the implicit function theorem, which, given the partial derivatives
$\Di_p A$ and $\Di_p f$, provides the derivative of the state $u$
w.r.t.\ $p$ as $\Di_p u = [\Di_u A]^{-1} (\Di_p f - \Di_p A)$.  This
--- and higher derivatives --- may be directly used in the
approximation of $u^*(p)$, as well as for a priori bounds for some
approximations.  These topics will not be pursued further here.

Furthermore assume that we are also given an iterative solver --- convergent
for all values of $p$ --- which generates successive iterates for $k=0,\ldots,$
\begin{equation} \label{eq:It}
   u^{(k+1)}(p) = S(p;u^{(k)}(p),R(p;u^{(k)}(p)),
  \quad\text{with}\; u^{(k)}(p) \rightarrow u^*(p),
\end{equation}
where $S$ is one cycle of the solver which may also depend on the iteration
counter $k$, $u^{(0)}$ is some starting vector, and $R(p;u^{(k)}(p))$
is the residuum of \refeq{eq:I}
\begin{equation} \label{eq:Res}
    R(u^{(k)}) := R(p;u^{(k)}(p)) := f(p) - A(p;u^{(k)}).
\end{equation}
Obviously, when the residuum vanishes --- $R(p;u^*(p))=0$ --- the
mapping $S$ has a fixed point $u^*(p) = S(p;u^*(p),0)$.

This mapping $S$ is the mathematical formalisation of the software
interface we will be assuming in order to derive a non-intrusive Galerkin
method, i.e.\ we will assume that the mapping $S$ is applied to its
inputs with one invocation of the ``solver''.

In the iteration in \refeq{eq:It} we may set $u^{(k+1)} = u^{(k)}+
\Delta u^{(k)}$ with
\begin{align} \label{eq:delta}
  \Delta u^{(k)} &:= S(p;u^{(k)},R(p;u^{(k)}))- u^{(k)}, \quad \text{and usually}\\
  \label{eq:It2}
  P (\Delta u^{(k)}) &= R(p;u^{(k)}),
\end{align}
so that in \refeq{eq:It}: $S(p;u^{(k)}) = u^{(k)} + P^{-1}(R(p;u^{(k)}))$,
where by slight abuse of notation we have shortened the list of
arguments.  Here $P$ is some preconditioner, which may depend on $p$,
the iteration counter $k$, and on the current iterate $u^{(k)}$; e.g.\
in \emph{Newton's method} $P = \Di_u A(p;u^{(k)})$.  In any case, we
assume that for all arguments the map $P$ is linear in $\Delta u$ and non-singular.
The iteration corresponding to a normal solve for a particular value
of $p$ then is given in Algorithm~\ref{alg:basic}.
 \begin{algorithm}
  \caption{Iteration of \refeq{eq:It}}
       \label{alg:basic}
   \begin{algorithmic}
    \State Start with some initial guess  $u^{(0)}$
    \State $k\gets 0$
    \While  {\emph{no convergence}} 
       \Statex\Comment{ \%comment: the global iteration loop\%} 
       \State Compute $\Delta u^{(k)}$ according to \refeq{eq:delta} or \refeq{eq:It2}
       \State $u^{(k+1)} \gets u^{(k)} + \Delta u^{(k)}$
       \State $k\gets k+1$
    \EndWhile 
  \end{algorithmic}
\end{algorithm}

We will assume additionally that the iteration converges at least
linearly, i.e.\ one has $\|\Delta u^{(k+1)}(p)\|_{\C{U}} \leq \vrho(p)\,
\|\Delta u^{(k)}(p)\|_{\C{U}}$, with $\vrho(p) < 1$.  For the convergence
analysis to follow later we will assume that the convergence factors
or Lipschitz constants $\vrho(p)$ are uniformly bounded for all values
of $p\in\C{P}$ by a constant strictly less than unity, i.e.\ $\vrho(p)
\leq \vrho^* < 1$.  Another way of saying this is that for all $u,v
\in \C{U}$ and $p\in\C{P}$ the iterator $S$ in \refeq{eq:It}
is uniformly Lipschitz continuous with Lipschitz constant $\vrho^* <
1$, i.e.\ a strict contraction:
\begin{equation} \label{eq:Lip}
    \|S(p;u(p),R(p;u(p))) - S(p;v(p),R(p;v(p)))\|_{\C{U}} \leq 
       \vrho^* \, \|u(p)-v(p)\|_{\C{U}}.
\end{equation}
One may recall from Banach's fixed point theorem that this provides us
with the a posteriori error bounds
\begin{equation}  \label{eq:contract}
 \|u^*(p) - u^{(k+1)}(p) \|_{\C{U}}  \leq \frac{\vrho^*}{1-\vrho^*}
 \|\Delta u^{(k)}(p) \|_{\C{U}},
\end{equation}
while the satisfaction of the equation may be gauged by
$\|R(p;u^{(k)})\|_{\C{U}}$.

%
%
%
%
%
%
%
%
%
%


%% file: chapters/galerkin.tex
%

\section{Galerkin approximation of parametric dependence}
   \label{S:galerkin}
The describe the dependence of $u$ on the parameters $p$ one would like to
approximate $u^*(p)$ in the following fashion:
\begin{equation} \label{eq:Gal-app}
   u^*(p) \approx u_{\C{I}}(p) = \sum_{\alpha \in \C{I}} u_{\alpha} \psi_{\alpha}(p),
\end{equation}
where $u_{\alpha} \in \C{U}$ are vector coefficients to be determined,
and $\psi_{\alpha}$ are some linearly independent functions, whose
linear combinations $\C{Q}_{\C{I}} := \spn\{\psi_{\alpha}\}_{\alpha
  \in \C{I}} \subset \D{R}^\C{P}$ form the Galerkin subspace of
parametric ``ansatz'' functions, and $\C{I}$ is some finite set of
(multi-)indices of cardinality $M:=|\C{I}|$.  Often the set $\C{I}$ has
no canonical order, but for the purpose of computation later we will
assume that some particular ordering has been chosen.

If we take the \emph{ansatz} \refeq{eq:Gal-app} and insert it into \refeq{eq:I},
the residuum \refeq{eq:Res} will usually not vanish for all $p$, as the finite
set of functions $\{\psi_{\alpha}\}_{\alpha \in \C{I}}$ can not match all possible
parametric variations of $u(p)$.

\subsection{The Galerkin equations for the residual} 
            \label{SS:galerkin-eq}
The \emph{Galerkin} method --- also called the method of
weighted residuals --- determines the unknown coefficients
$u_{\alpha}$ in \refeq{eq:Gal-app} by requiring that for
all $\beta \in \C{I}$
\begin{equation} \label{eq:Gal-cond}
\vek{G}_{\C{Q}}(\vphi_\beta(\cdot) R(\cdot;u_{\C{I}})) = 0,
\end{equation}
where  $\{\vphi_{\beta}\}_{\beta \in \C{I}}$
is a set of linearly independent functions of $p$.  The residuum
$R(p;u_{\C{I}}(p))$ in the argument of the linear Galerkin projector
$\vek{G}_{\C{Q}}$ is a parametric function, and such functions may be
represented by a sum $R(p;u_{\C{I}}(p))=\sum_n \phi_n(p) v_n$ with $\phi_n \in
\D{R}^{\C{P}}$.  Hence the projector may be defined by requiring that for
scalar functions $\psi, \phi \in \C{Q} \subseteq \D{R}^{\C{P}}$ and
a vector $v \in \C{U}$ one has
\begin{equation} \label{eq:Gal-pro}
\vek{G}_{\C{Q}}(\phi(\cdot)\psi(\cdot)\, v) = \ip{\phi}{\psi}_{\C{Q}}\, v,
\end{equation}
where $\ip{\cdot}{\cdot}_{\C{Q}}$ is some duality pairing or inner
product on a subspace $\C{Q}$ of the scalar functions, and from this
$\vek{G}_{\C{Q}}$ can be extended by linearity:
$$\vek{G}_{\C{Q}}(\vphi_\beta\, R(\cdot;u_{\C{I}})) =
\vek{G}_{\C{Q}}(\vphi_\beta\, \sum_n \phi_n\, v_n) = \sum_n
\vek{G}_{\C{Q}}(\vphi_\beta\, \phi_n\, v_n) = \sum_n
\ip{\vphi_\beta}{\phi_n}_{\C{Q}}\, v_n.$$ 
It is easy to see that this definition is independent of the
particular representation of the parametric function.

In case $\C{P}$ is a measure space with measure $\mu$, then that
pairing often is $\ip{\phi}{\psi}_{\C{Q}} = \int_{\C{P}} \phi(p)
\psi(p) \, \mu(\di p)$, and if $\mu(\C{P}) = 1$, such that $\C{P}$
may be considered as a probability space with expectation operator
$\EXP{\phi} = \int_{\C{P}} \phi(p)\, \mu(\di p)$, then
$\ip{\phi}{\psi}_{\C{Q}} = \EXP{\phi \psi}$.  Observe that a sum
like $\ip{\phi}{\psi}_{\C{Q}} = \sum_j w_j \phi(p_j) \psi(p_j)$ with
positive weights $w_j$ is a special form of such an integral.  A bit
more general would be to allow $\ip{\phi}{\psi}_{\C{Q}}
=\iint_{\C{P}\times\C{P}} \varkappa(p,q) \phi(p) \psi(q) \, \mu(\di p)
\mu(\di q)$, where $\varkappa$ is a symmetric positive definite
kernel.  What is important for what is to follow, and what we want to
assume from now on, is that the pairing is given by some integral, and
we will assume the form $\int_{\C{P}} \phi(p) \psi(p) \, \mu(\di p)$
for the sake of simplicity.

The set $\{\psi_{\alpha}\}_{\alpha \in \C{I}}$ determines the
Galerkin subspace $\C{Q}_{\C{I}} := \spn\{\psi_{\alpha}\}_{\alpha \in
  \C{I}} \subseteq \C{Q}$, which is responsible for the approximation
properties, whereas the set $\{\vphi_{\beta}\}_{\beta \in \C{I}}$
determines the projection onto that subspace which is important for
the stability of the procedure, as the projection is orthogonal to
$\hat{\C{Q}}_{\C{I}} := \spn\{\vphi_{\beta}\}_{\beta \in \C{I}}$.
Often one takes $\vphi_\beta = \psi_\beta$ and hence
$\hat{\C{Q}}_{\C{I}}= \C{Q}_{\C{I}}$, and this is then commonly called
the \emph{Bubnov-Galerkin} method, whereas in the general case
$\hat{\C{Q}}_{\C{I}} \neq \C{Q}_{\C{I}}$ one speaks of the
\emph{Petrov-Galerkin} method.

Explicitly writing down \refeq{eq:Gal-cond}, one obtains for all $\beta$
\begin{equation} \label{eq:Gal-cond-expl}
\vek{G}_{\C{Q}} (\vphi_\beta(\cdot) (f(p) - A(p;\sum_{\alpha \in \C{I}} 
                 u_{\alpha} \psi_{\alpha}(p)))) = 0.
\end{equation}
It is important to recognise that \refeq{eq:Gal-cond-expl} is a ---
usually coupled --- system of equations for the unknown vectors
$u_\alpha$ of size $M \times N$, as $M = \dim \C{Q}_{\C{I}}$ and $N =
\dim \C{U}$.  These equations look sufficiently different from
\refeq{eq:I}, so that the common wisdom is that the solution of
\refeq{eq:Gal-cond-expl} requires new software and new methods, and
that the solver \refeq{eq:It} is of no use here.  As a change or
re-write of the existing software seems to be necessary, the resulting
methods are often labelled ``intrusive''.

As a remark, observe that if one chooses $\vphi_\beta(p) =
\delta_\beta(p) = \delta(p - p_\beta)$ --- the delta-``function''
associated to the duality pairing $\ip{\cdot}{\cdot}_{\C{Q}}$ (i.e.\
$\ip{\delta_\beta}{\phi}_{\C{Q}} = \phi(p_\beta)$) --- where the $p_\beta$
are distinct points in $\C{P}$ in \refeq{eq:Gal-it}, this becomes for
all $\beta$:
\begin{multline} \label{eq:coll-1}
  0=\vek{G}_{\C{Q}}(\delta_\beta R(\cdot;u_{\C{I}})) = R(p_\beta;u_{\C{I}}(p_\beta)) 
   = \\ f(p_\beta) - A(p_\beta;\sum_{\alpha \in \C{I}} u_{\alpha} \psi_{\alpha}(p_\beta))
    = f(p_\beta) - A(p_\beta; u_{\beta}),
\end{multline}
where the last of these equalities holds only in case the basis
$\{\psi_\alpha\}$ satisfies the \emph{Kronecker}-$\delta$ property
$\psi_\alpha(p_\beta) = \delta_{\alpha, \beta}$, as then $u_\beta =
u_{\C{I}}(p_\beta)$.  In this latter case these are $M$ uncoupled
equations each of size $N$, and they have for each $p_\beta$ the form
\refeq{eq:I} --- we have recovered the \emph{collocation} method which
independently for each $p_\beta$ computes $u_\beta$, using the solver
\refeq{eq:It}.  Such a method then obviously is non-intrusive, as the
original software may be used.  Thus this is often the method of
choice, as often there is considerable investment in the software
which performs \refeq{eq:It}, which one would like to re-use.
Unfortunately this choice is very rigid as regards the subspace
$\C{Q}_{\C{I}}$ and the projection orthogonal to
$\hat{\C{Q}}_{\C{I}}$.

We believe that this is a \emph{false} alternative, and that the
distinction is not between \emph{intrusive} or \emph{non-intrusive},
but between \emph{coupled} or \emph{uncoupled}.  Furthermore, and more
importantly, we want to show that also in the more general case of a
coupled system like in \refeq{eq:Gal-cond-expl} the original solver
\refeq{eq:It} may be put to good use.  This will be achieved by making
\refeq{eq:It} the starting point, instead of \refeq{eq:I} or
\refeq{eq:Res}.  Such coupled iterations also arise for example from
multi-physics problem, and there these coupled iterations can also be
solved by what is called a \emph{partitioned} approach, see e.g.\
\cite{matthiesStf03}, which is the equivalent of non-intrusive here.
Quite a number of different variants of global partitioned iterations
are possible \cite{matthiesStf03}, we only look at some of the
simplest variants, as the point is here only to dispel the myth about
intrusiveness.

\subsection{The fixed-point Galerkin equations} \label{SS:galerkin-fixed}
Whatever the starting point, we would still like to achieve the same
result.  So before continuing, let us show
\begin{prop} \label{prop:project} 
  Projecting the fixed point equation attached to the iteration \refeq{eq:It},
  namely $u_{\C{I}}=u_{\C{I}} + P^{-1}(R(u_{\C{I}}))$, is equivalent to
  projecting the preconditioned residual $P^{-1}(R(u_{\C{I}}))$, that means for
  all $\beta \in \C{I}$
 \begin{equation}  \label{eq:prop-1}
  \vek{G}_{\C{Q}} (\vphi_\beta(\cdot) \, P^{-1}(\cdot)(R(\cdot;u_{\C{I}}(\cdot)))) = 0.
 \end{equation}
  Moreover, if the preconditioner $P$ in \refeq{eq:It2} does not depend on $p$ nor $u$,
  then it is equivalent to projecting the residual $R(u_{\C{I}})$ from 
  \refeq{eq:Gal-cond}, that means for all $\beta \in \C{I}$
 \begin{equation}  \label{eq:prop-2}
  \vek{G}_{\C{Q}} (\vphi_\beta(\cdot) \, R(\cdot;u_{\C{I}}(\cdot))) = 0.
 \end{equation}
\end{prop}
\begin{proof} 
The \refeq{eq:prop-1} follows simply from linearity of $\vek{G}_{\C{Q}}$.  Furthermore,
in case $P$ does not depend on $p$ nor $u$, for \refeq{eq:prop-2} we have from
\refeq{eq:prop-1} for any $\beta \in \C{I}$
\begin{equation}  \label{eq:prop}
0 = \vek{G}_{\C{Q}} (\vphi_\beta(\cdot) \, P^{-1}(R(\cdot;u_{\C{I}}(\cdot)))) = 
  P^{-1} \vek{G}_{\C{Q}} (\vphi_\beta\,R(u_{\C{I}})) \quad
\Leftrightarrow \quad 0 = \vek{G}_{\C{Q}} (\vphi_\beta\,R(u_{\C{I}})),
\end{equation}
on noting that for any linear map $L$ one has $\vek{G}_{\C{Q}}(\vphi(\cdot)\,
L (\phi(\cdot)\, v)) = \ip{\vphi}{\phi}_{\C{Q}}\, L v = L\,
\vek{G}_{\C{Q}}(\vphi(\cdot)\, \phi(\cdot)\, v)$, and by observing that 
$P^{-1}$ is non-singular.
\end{proof}

This means that instead of the residual \refeq{eq:Res} we may just as
well project the iteration \refeq{eq:It2}: with the abbreviation 
$R^{(k)}(\cdot):= R(\cdot;u^{(k)}(\cdot))$ we have for all $\beta \in \C{I}$
\begin{equation} \label{eq:Gal-it}
\vek{G}_{\C{Q}}(\vphi_\beta(\cdot)\, u^{(k+1)}) = 
    \vek{G}_{\C{Q}}(\vphi_\beta \, (u^{(k)} + \Delta u^{(k)})) =
    \vek{G}_{\C{Q}}(\vphi_\beta \, (u^{(k)} + P^{-1}R^{(k)})).
\end{equation}
Expanding $u^{(k)}(p) = \sum_\alpha u^{(k)}_\alpha \psi_\alpha(p)$ in 
\refeq{eq:Gal-it}, that becomes a coupled iteration equation for the $u_\alpha$:
\begin{equation} \label{eq:Gal-it-alpha}  \forall \beta:\;
\vek{G}_{\C{Q}}(\vphi_\beta(\cdot)\, \sum_\alpha u^{(k+1)}_\alpha \psi_\alpha(\cdot)) =
    \vek{G}_{\C{Q}}(\vphi_\beta(\cdot)\, (\sum_\alpha u^{(k)}_\alpha \psi_\alpha(\cdot) 
    + P^{-1}R^{(k)}(\cdot))),
\end{equation}
which may now be written as
\begin{equation} \label{eq:Gal-it-1}  \forall \beta: \;
\sum_\alpha \vek{M}_{\beta,\alpha} u^{(k+1)}_\alpha = 
\sum_\alpha \vek{M}_{\beta,\alpha} u^{(k)}_\alpha +
    \vek{G}_{\C{Q}}(\vphi_\beta\, P^{-1}R^{(k)}),
\end{equation}
where $\vek{M}_{\beta,\alpha} :=
\ip{\vphi_\beta}{\psi_\alpha}_{\C{Q}}$.  If the coefficients $u^{(k)}_\alpha
\in \C{U}$ are arranged column-wise in a $N \times M$ matrix $\mat{u}^{(k)}
= [\dots, u^{(k)}_\alpha, \dots] \in \C{U}^{\C{I}}$ and similarly
$\mat{G}_{\C{Q}}(P^{-1}R^{(k)}) = [\dots,
\vek{G}_{\C{Q}}(\vphi_\alpha \, P^{-1}R^{(k)}), \dots]$, and
the $\vek{M}_{\beta,\alpha}$ are viewed as entries of a $M \times M$
matrix $\mat{M} \in \D{R}^{\C{I}\times\C{I}}$, \refeq{eq:Gal-it-1} may
be compactly written as
\begin{align} \label{eq:Gal-it-1aa}
\mat{u}^{(k+1)} \mat{M}^T &= \mat{u}^{(k)} \mat{M}^T + 
 \mat{G}_{\C{Q}}(P^{-1}R^{(k)}), \quad\text{  or as}\\  \label{eq:Gal-it-1c}
\mat{u}^{(k+1)}  &= \mat{u}^{(k)}  + \vek{\Delta}_{\C{Q}}(\mat{u}^{(k)})
                =: \mat{S}_{\C{Q}}(\mat{u}^{(k)}),
\end{align}
where we have defined two new functions
$\vek{\Delta}_{\C{Q}}(\mat{u}^{(k)}) := [\mat{G}_{\C{Q}}(P^{-1}R^{(k)})] 
\mat{M}^{-T}$, and $\mat{S}_{\C{Q}}(\mat{u}^{(k)}) = \mat{u}^{(k)} +
\vek{\Delta}_{\C{Q}}(\mat{u}^{(k)})$, which will be needed later for
the convergence analysis in \refSS{coupled-conv}.  

It is apparent that the computation will be much simplified if the
ansatz-functions $\{\psi_{\alpha}\}_{\alpha \in \C{I}}$ and the
test-functions for the projection $\{\vphi_{\beta}\}_{\beta \in \C{I}}$
are chosen bi-orthogonal, i.e.\ if one has for all $\alpha, \beta \in
\C{I}$ that $\vek{M}_{\beta, \alpha} = \delta_{\beta, \alpha}$, i.e.\
$\mat{M} = \mat{I}$, which shall be assumed from now on.  Hence now
\begin{equation}  \label{eq:delta-BJ}
  \vek{\Delta}_{\C{Q}}(\mat{u}^{(k)}) =
  \mat{G}_{\C{Q}}(P^{-1}R^{(k)}) = 
   [\dots, \vek{G}_{\C{Q}}(\vphi_\alpha\,P^{-1}R^{(k)}), \dots].
\end{equation}

\refeq{eq:Gal-it-1c} is already a possible way of performing
the iteration.  The practical, non-intrusive, computation of the terms
in \refeq{eq:Gal-it-1c} still has to be considered, but we may
formulate the corresponding Algorithm~\ref{alg:bJ} and investigate its
convergence beforehand.  The reader who is only interested in the
computational description of the non-intrusive algorithm may jump
directly to \refS{resid}.
\begin{algorithm}
  \caption{Block Jacobi iteration of \refeq{eq:Gal-it-1c}}
    \label{alg:bJ}
   \begin{algorithmic}
     \State Start with some initial guess $\mat{u}^{(0)}$
     \State $k\gets 0$ 
     \While {\emph{no convergence}} 
        \Statex\Comment{ \%comment: the global iteration loop\%} 
        \State Compute $\vek{\Delta}_{\C{Q}}(\mat{u}^{(k)})$ according 
                to \refeq{eq:delta-BJ} 
        \State $\mat{u}^{(k+1)} \gets \mat{u}^{(k)}+\vek{\Delta}_{\C{Q}}(\mat{u}^{(k)}) 
            \quad [ =  \mat{S}_{\C{Q}}(\mat{u}^{(k)}) ]$ 
        \State $k\gets k+1$
     \EndWhile 
  \end{algorithmic}
\end{algorithm}

Although the underlying iteration \refeq{eq:It} in
Algorithm~\ref{alg:basic} may be of any kind --- e.g.\ Newton's method
--- when one views \refeq{eq:Gal-it-1c} with regard to the block
structure imposed by the $\mat{u} = [\dots, u_\beta, \dots]$,
Algorithm~\ref{alg:bJ} is a --- maybe nonlinear --- \emph{block
  Jacobi} iteration.

\subsection{Convergence of coupled iterations} \label{SS:coupled-conv}
Here we want to show that the map $\mat{S}_{\C{Q}}$ in
\refeq{eq:Gal-it-1c} satisfies a Lipschitz condition with the same
constant as in \refeq{eq:Lip}.  This will need some more theoretical
considerations.  For the sake of simplicity we will assume that
$\ip{\cdot}{\cdot}_{\C{Q}}$ is actually an inner product on the
Hilbert space $\C{Q} \subseteq \D{R}^{\C{P}}$, such that
$\C{Q}_{\C{I}} \subseteq \C{Q}$.  The contraction condition for
$\mat{S}_{\C{Q}}$ with contraction factor (Lipschitz constant) less or
equal to $\vrho^*$ will only hold if the Galerkin projection is
\emph{orthogonal}, i.e.\ we have to take $\vphi_\alpha = \psi_\alpha$,
which means $\hat{\C{Q}}_{\C{I}} = \C{Q}_{\C{I}}$.  Our previous
assumption that $\mat{M} = \mat{I}$ --- which is now the Gram matrix
of the basis $\{\psi_\alpha\}_{\alpha \in \C{I}}$ --- now means that
this basis is actually \emph{orthonormal}.

Parametric elements like $\C{P} \ni p\mapsto u(p) \in \C{U}$ are
formally in the Hilbert tensor product space of sums like $\sum_n
\phi_n(p) v_n =: \sum \phi_n \otimes v_n \in \C{Q} \otimes \C{U}$,
with the inner product of two elementary tensors $\phi_j \otimes w_j
\in \C{Q} \otimes \C{U}, (j=1,2)$, defined by $\ip{\phi_1 \otimes
  w_1}{\phi_2 \otimes w_2}_{\C{Q} \otimes \C{U}} :=
\ip{\phi_1}{\phi_2}_{\C{Q}} \ip{w_1}{w_2}_{\C{U}}$, and then extended
by bi-linearity.  In the space $\C{Q}$, the subspace $\C{Q}_{\C{I}}$,
which is finite-dimensional and hence closed,
leads to the orthogonal direct sum decomposition $\C{Q} =
\C{Q}_{\C{I}} \oplus \C{Q}_{\C{I}}^\perp$, and hence to the orthogonal
direct sum decomposition $\C{Q} \otimes \C{U} = (\C{Q}_{\C{I}} \otimes
\C{U}) \oplus (\C{Q}_{\C{I}}^\perp \otimes \C{U})$.

The mapping $J: \C{U}^{\C{I}} \ni \mat{u} = [\dots, u_\alpha, \dots]
\mapsto \sum_\alpha \psi_\alpha(\cdot) u_\alpha \in \C{Q}_{\C{I}}
\otimes \C{U} \subseteq \C{Q}\otimes\C{U}$ is by design bijective
onto $\C{Q}_{\C{I}} \otimes \C{U} $ and
may thus be used to induce a norm and inner product on $\C{U}^{\C{I}}$
via $\|\mat{u}\|_{\C{U}^{\C{I}}}^2 := \|J \mat{u} \|_{\C{Q}\otimes\C{U}}^2
= \|\sum_\alpha\psi_\alpha(\cdot) u_\alpha \|_{\C{Q}\otimes\C{U}}^2 = 
\sum_\alpha\|u_\alpha\|_{\C{U}}^2$, making it a unitary map, hence $\| J\| = 1$.
When viewed as a mapping into the larger space $\C{Q}\otimes\C{U}$,
were it is extended by slight abuse of notation by the inclusion,
it remains an isometry.

\begin{lem}  \label{lem:G_Q-T}
The maps $\mat{G}_{\C{Q}}: \C{Q}\otimes\C{U} \to \C{U}^{\C{I}}$ and 
$J$ are adjoints of each other, $ \mat{G}_{\C{Q}}^*
= J$, and $\mat{G}_{\C{Q}}$ is non-expansive,
$\|\mat{G}_{\C{Q}}\| = \|\mat{G}_{\C{Q}}^*\| = 1$.
\end{lem}
\begin{proof}
 For all $\mat{v} \in \C{U}^{\C{I}}$ and $\phi \otimes w \in \C{Q}\otimes\C{U}:$
\begin{multline}  \label{eq:G_Q-trans-1}
   \ip{\mat{G}_{\C{Q}}(\phi \otimes w)}{\mat{v}}_{\C{U}^{\C{I}}} = 
   \ip{[\ldots, \ip{\psi_\alpha}{\phi}_{\C{Q}}\, w, \ldots]}{[\ldots,
  v^\alpha, \ldots]}_{\C{U}^{\C{I}}} = \\
  \sum_\alpha\ip{\psi_\alpha}{\phi}_{\C{Q}}  \ip{w}{v^\alpha}_{\C{U}}  =
  \ip{\phi \otimes w}{\sum_\alpha \psi_\alpha \otimes v^\alpha}_{\C{Q}\otimes\C{U}}
  = \ip{\phi \otimes w}{J \mat{v}}_{\C{Q}\otimes\C{U}},
\end{multline}
and hence $\mat{G}_{\C{Q}}^* = J$.  But $J$ is an isometry, so that one has
$\|\mat{G}_{\C{Q}}^*\| = \| J\|  = 1$.
As $\|\mat{G}_{\C{Q}}\| = \|\mat{G}_{\C{Q}}^*\|$, we are finished.
\end{proof}

With the observation that
\begin{multline}  \label{eq:G_Q-extens}
  \mat{G}_{\C{Q}}(S(\cdot;u^{(k)}(\cdot),R^{(k)}(\cdot))= 
  \mat{G}_{\C{Q}}(u^{(k)}(\cdot) + P^{-1}R^{(k)}(\cdot))=\\
  \mat{u}^{(k)} + \mat{G}_{\C{Q}}(P^{-1}R^{(k)}) = \mat{u}^{(k)} +
  \vek{\Delta}_{\C{Q}}(\mat{u}^{(k)}) = \mat{S}_{\C{Q}}(\mat{u}^{(k)}),
\end{multline}
the map $\mat{S}_{\C{Q}}: \C{U}^{\C{I}} \to \C{U}^{\C{I}}$ in
\refeq{eq:Gal-it-1c} may be factored in the following way:
\begin{align}
  \label{eq:fact-map-dia}
  \mat{S}_{\C{Q}} &: \C{U}^{\C{I}} \overset{J}{\to} \C{Q}\otimes\C{U}
   \overset{\tilde{S}}{\to} \C{Q}\otimes\C{U}
   \overset{\mat{G}_{\C{Q}}}{\to} \C{U}^{\C{I}},\\
     \label{eq:fact-map-S}
  \mat{S}_{\C{Q}} &= \mat{G}_{\C{Q}} \circ \tilde{S} \circ J
  = \mat{G}_{\C{Q}} \circ \tilde{S} \circ \mat{G}_{\C{Q}}^*,
\end{align}
where $\tilde{S}$ is defined via the solver map $S$ in \refeq{eq:It} by
\begin{equation}   \label{eq:s-til}
  \tilde{S}: \C{Q}\otimes\C{U} \ni u(\cdot) \mapsto S(\cdot; u(\cdot), 
  R(\cdot, u(\cdot))) \in \C{Q}\otimes\C{U}.
\end{equation}
For this mapping we have the following result:
\begin{prop}  \label{prop:contr-s-til}
  In \refeq{eq:s-til}, the map denoted $\tilde{S}$ has the same Lipschitz
  constant $\vrho^*$ as the map $S$ in \refeq{eq:It}, cf.\
  \refeq{eq:Lip}; i.e.\ $\tilde{S}$  is a contraction with contraction
  factor $\vrho^* < 1$. 
\end{prop}
\begin{proof}
  We now use the assumption that the inner product on $\C{Q}$ is given
  by an integral,
  $\ip{\vphi}{\phi}_{\C{Q}} = \int_{\C{P}} \vphi(p) \phi(p)\, \mu(\di p)$.
  In that case $\C{Q} = L_2(\C{P},\mu;\D{R})$, and the Hilbert tensor
  product $\C{Q}\otimes\C{U}$ is isometrically isomorphic to
  $L_2(\C{P},\mu;\C{U})$.  Hence with \refeq{eq:Lip} for all
  $u(\cdot), v(\cdot) \in L_2(\C{P},\mu;\C{U})$
\begin{multline*}
\| \tilde{S}(u(\cdot)) - \tilde{S}(v(\cdot))\|^2_{L_2(\C{P},\mu;\C{U})} =\\
\int_{\C{P}}\|S(p;u(p),R(p;u(p))) - S(p;v(p),R(p;v(p)))\|^2_{\C{U}}\,\mu(\di p)\\ 
   \leq (\vrho^*)^2  \int_{\C{P}} \|u(p) - v(p)\|^2_{\C{U}}\, \mu(\di p) =
  (\vrho^*)^2 \| u(\cdot) - v(\cdot)\|^2_{L_2(\C{P},\mu;\C{U})},
\end{multline*}
and the proof is concluded by taking square roots.
\end{proof}
This immediately leads to
\begin{coro} \label{coro:Lip-s} 
  The map $\mat{S}_{\C{Q}}$ from \refeq{eq:Gal-it-1c} is a contraction with
  contraction factor $\vrho^* < 1$ (see \refeq{eq:Lip}):
\begin{equation} \label{eq:Lip-BJ}  \forall \mat{u}, \mat{v} \in
  \C{U}^{\C{I}}: \quad
    \|\mat{S}_{\C{Q}}(\mat{u}) - \mat{S}_{\C{Q}}(\mat{v})\|_{\C{U}^{\C{I}}} \leq 
     \vrho^* \, \|\mat{u}-\mat{v}\|_{\C{U}^{\C{I}}},
\end{equation}
and hence the Galerkin equations have a unique solution $\mat{u}^* \in
\C{U}^{\C{I}}$.
\end{coro}
\begin{proof}
This follows from the decomposition \refeq{eq:fact-map-S}, Lemma~\ref{lem:G_Q-T},
and Proposition~\ref{prop:contr-s-til}, as $\|\mat{S}_{\C{Q}}\| =
\|\mat{G}_{\C{Q}}\circ \tilde{S} \circ \mat{G}_{\C{Q}}^*\| \leq
\|\mat{G}_{\C{Q}}\|\,  \|\tilde{S} \|\, \| \mat{G}_{\C{Q}}^*\| \leq
\vrho^*$, and Banach's contraction mapping theorem.
\end{proof}

Now we may state the main result about the convergence of the simple
\emph{block-Jacobi} Algorithm~\ref{alg:bJ} for the coupled Galerkin system:
\begin{thm} \label{thm:contract-s-p} 
  As the map $\mat{S}_{\C{Q}}$ from \refeq{eq:Gal-it-1c} has Lipschitz
  constant $\vrho^* < 1$, and is thus a contraction with the same
  factor as the solver $S$ in \refeq{eq:It}, the
  Algorithm~\ref{alg:bJ} converges to the unique solution $\mat{u}^*
  \in \C{U}^{\C{I}}$ with the same linear speed of convergence as
  Algorithm~\ref{alg:basic}.  Additionally, we have the a posteriori
  error estimate (see \refeq{eq:contract})
\begin{equation}  \label{eq:contract-BJ}
   \|\mat{u}^* - \mat{u}^{(k+1)} \|_{\C{U}^{\C{I}}}  \leq \frac{\vrho^*}{1-\vrho^*} \;
              \|\vek{\Delta}_{\C{Q}}(\mat{u}^{(k)}) \|_{\C{U}^{\C{I}}}.
\end{equation}
The satisfaction of the parametric equation may be gauged by
$\|R^{(k)} \|_{\C{Q}\otimes\C{U}} = \|R(\cdot;u^{(k)}) \|_{\C{Q}\otimes\C{U}}$.
\end{thm}
\begin{proof}
Everything simply follows from Corollary~\ref{coro:Lip-s}, Banach's
contraction mapping theorem, and the fact that $R^{(k)}(\cdot)$
is the residuum at iteration $k$ before any preconditioning or projection.
\end{proof}
Observe that this only holds for the linear convergence speed; in case 
Algorithm~\ref{alg:basic} has super-linear convergence, this can not
be necessarily matched by Algorithm~\ref{alg:bJ}, for this more
sophisticated algorithms for the coupled equations are necessary, 
see e.g.\ \cite{matthiesStf03}.

%
%
%
%
%
%
%
%
%
%
%


%% file: chapters/residual.tex
%

\section{Non-intrusive residual} \label{S:resid}
Here we want to look in more detail at the actual computation of the
right hand side of \refeq{eq:Gal-it-1c}, in the form \refeq{eq:delta-BJ}.
One may observe that the term $\vek{G}_{\C{Q}}(\vphi_\alpha\, P^{-1}
R^{(k)})$ in \refeq{eq:delta-BJ} is the Galerkin projection
of the preconditioned residual for that iteration.  
Let us recall that the Galerkin projector was defined by
\refeq{eq:Gal-pro} as
$\vek{G}_{\C{Q}}(\vphi(\cdot) \,v\, \phi(\cdot)) = \ip{\vphi}{\phi}_{\C{Q}} \, v$.

\subsection{Analytic computation}   \label{SS:AnalytComp}
In some cases \cite{boulder:2011, HgmWN12}, notably when the preconditioner $P$
does not depend on $p$ nor $u$, or when the operator $A$ is linear 
or polynomial in $u$, and linear in the parameters $p$, 
it may be possible to actually represent
$P^{-1} R^{(k)}$, not just in principle, but actually
\emph{non-intrusively} through the use of the solver $S$ in \refeq{eq:It} as
 \begin{equation}   \label{eq:residAnalytic}
   P^{-1} R^{(k)}(p) = \sum_n \rho_n(p)\, v_n
    = \sum_{n,\beta} \rho_{n,\beta}\, \psi_\beta(p)\, v_n,
 \end{equation}
where $\rho_{n,\beta} = \ip{\psi_\beta}{\rho_n}_{\C{Q}}$,  --- remembering that
we chose $\vphi_\alpha = \psi_\alpha$ orthonormal.
The Galerkin projection of this then gives
\begin{multline}  \label{eq:GPresAnalytic2}
  \vek{G}_{\C{Q}}(\psi_\alpha\, P^{-1}R^{(k)}) =
  \vek{G}_{\C{Q}}(\psi_\alpha \sum_{n,\beta} \rho_{n,\beta}
  \psi_\beta\, v_n) = \\ \sum_{n,\beta} \rho_{n,\beta} \,
  \vek{G}_{\C{Q}}(\psi_\alpha \psi_\beta\, v_n) = \sum_{n,\beta} \rho_{n,\beta} 
  \,\ip{\psi_\alpha}{\psi_\beta}_{\C{Q}} \, v_n =  \sum_n \rho_{n,\alpha}\, v_n,
\end{multline}
using the linearity of $\vek{G}_{\C{Q}}$ and the orthonormality of the
basis $\{ \psi_\alpha\}_{\alpha\in\C{I}}$.  This means that for the
right hand side of  \refeq{eq:Gal-it-1c} in the form \refeq{eq:delta-BJ},
given the representation \refeq{eq:residAnalytic}, each term may be
computed through simple linear algebra operations \refeq{eq:GPresAnalytic2}.
This expression may be directly used in the \emph{block-Jacobi}
Algorithm~\ref{alg:bJ} for $\vek{\Delta}_{\C{Q}}(\mat{u}^{(k)})$
in the form \refeq{eq:delta-BJ}, and the description of the algorithm
is complete.  Let us remark finally that if the solver actually returns 
$S(p;u^{(k)}(p),R^{(k)}(p))$ instead of the increment
$P^{-1} R^{(k)}(p)$, Algorithm~\ref{alg:bJ} is easily adapted
by computing completely analogously $\mat{S}_{\C{Q}}(\mat{u}^{(k)})$.

\subsection{Numerical integration}   \label{SS:NumerIntegr}
The following idea to obtain a non-intrusive computation of the right
hand side of  \refeq{eq:Gal-it-1c} in the form \refeq{eq:delta-BJ}, is
more general, but involves a further approximation, namely numerical
integration.

Remembering that it was assumed that the duality pairing on the scalar
functions is given by an integral with measure $\mu$,
\begin{equation} \label{eq:pairing} 
  \ip{\vphi}{\phi}_{\C{Q}} = \int_{\C{P}} \vphi(p) \phi(p) \, \mu(\di p),
\end{equation}
we now assume that this integral has some approximate numerical
quadrature formula 
\begin{equation} \label{eq:quad} 
  \int_{\C{P}} \phi(p) \, \mu(\di p) \approx \sum_{z=1}^Z w_z \phi(p_z),
\end{equation}
where the integrand is evaluated at the quadrature points $p_z$ and
the $w_z$ are appropriate weights.  

With this approximation the term $\vek{G}_{\C{Q}}(\psi_\beta\,
P^{-1}R^{(k)})$ in \refeq{eq:Gal-it-1} becomes practically
computable without any further assumptions on the operator $A$, giving
\begin{align} \label{eq:q-resi} 
  \vek{G}_{\C{Q}}(\psi_ \beta\, P^{-1}R^{(k)}) \approx 
  \Delta_{Z, \beta} u^{(k)} :=&
     \sum_z w_z \psi_\beta(p_z)\, \Delta u^{(k)}_z, \quad  \text{where}\\
     \label{eq:q-resi-2}\Delta u^{(k)}_z :=
   P^{-1}(p_z) R(p_z;u^{(k)}(p_z)) =& P^{-1}(p_z)
   \left(f(p_z) - A(p_z; u^{(k)}(p_z))\right), \quad \text{or}\\
   =& S(p_z;u^{(k)}(p_z),R(p_z;u^{(k)}(p_z))) - u^{(k)}(p_z)
   \label{eq:q-resi-3}
\end{align}
is the preconditioned residuum evaluated at $p_z$, and $u^{(k)}(p_z)=
\sum_\alpha u_\alpha^{(k)} \psi_\alpha(p_z)$.  This is indeed the only
interface needed to the original equation, something which can be
easily evaluated \emph{non-intrusively} as the iteration increment
$\Delta u^{(k)}_z$ in \refeq{eq:q-resi-2} in case the current state is
given as $u^{(k)}(p_z)$ for the parameter value $p_z$.  An alternative
form is given in \refeq{eq:q-resi-3}, which is one iteration of the
solver, starting at $u^{(k)}(p_z)$ for the parameter $p_z$.  This
variant is for the case when the solver actually returns 
$S(p; u^{(k)}(p), R^{(k)}(p))$ instead of the increment
$P^{-1} R^{(k)}(p)$.

\subsection{Non-intrusive iteration}   \label{SS:NonIntrIter}
The term in \refeq{eq:Gal-it-1c} in the form of \refeq{eq:delta-BJ}
has to be computed \emph{non-intrusively}.  Following
\refSS{NumerIntegr} about numerical integration of the terms --- if
applicable, the analytic counterpart from \refSS{AnalytComp} can be
easily substituted --- we formulate the approximation of
\begin{equation}    \label{eq:dqnumint}
  \vek{\Delta}_{\C{Q}}(\mat{u}^{(k)}) =
   [\dots, \vek{G}_{\C{Q}}(\vphi_\alpha\,P^{-1}R^{(k)}), \dots]
  \approx \vek{\Delta}_{Z}(\mat{u}^{(k)}) = 
   [\dots, \Delta_{Z, \alpha} u^{(k)}, \dots]
\end{equation}
in Algorithm~\ref{alg:bJ} from \refeq{eq:delta-BJ} in
Algorithm~\ref{alg:NI}, using \refeq{eq:q-resi} and \refeq{eq:q-resi-2}: 
\begin{algorithm}
  \caption{Non-intrusive computation of \refeq{eq:delta-BJ} in the 
    form of \refeq{eq:dqnumint}}
    \label{alg:NI}
   \begin{algorithmic}
     \For {$\alpha \in \C{I}$}  
          \State $\Delta_{Z, \alpha} u^{(k)}\gets 0 $ 
     \EndFor
        \Statex\Comment{ \%comment: the loop over integration points\%} 
     \For {$z\gets 1,\ldots,Z$} 
        \State Compute $\Delta u^{(k)}_z$ from \refeq{eq:q-resi-2}
        \State $r_z \gets w_z\, \Delta u^{(k)}_z$  
        \For {$\alpha \in \C{I}$}
          \State $\Delta_{Z, \alpha} u^{(k)}\gets \Delta_{Z, \alpha} u^{(k)} +
                 \psi_\alpha(p_z) \, r_z$ 
        \EndFor
     \EndFor
  \end{algorithmic}
\end{algorithm}

The result of this algorithm is $\vek{\Delta}_{Z}(\mat{u}^{(k)})$, the
approximation of $\vek{\Delta}_{\C{Q}}(\mat{u}^{(k)})$ by numerical
integration.  With Algorithm~\ref{alg:NI} it is now possible to
formulate a non-intrusive version of Algorithm~\ref{alg:bJ}, the
block-Jacobi iteration, in Algorithm~\ref{alg:BJNI}.
 \begin{algorithm}
   \caption{Non-Intrusive block Jacobi iteration of \refeq{eq:Gal-it-1c}}
   \label{alg:BJNI}
   \begin{algorithmic}
    \State Start with some initial guess 
         $\Mtil{u}^{(0)}=[\dots,\tilde{u}^{(0)}_\alpha,\dots]$
    \State $k \gets 0$
    \While  {\emph{no convergence}} 
    \Statex\Comment{ \%comment: the global iteration loop\%} 
       \State Compute $\vek{\Delta}_{Z}(\Mtil{u}^{(k)}) =
                   [\dots,\Delta_{Z, \alpha} \tilde{u}^{(k)},\dots]$ 
              according to  Algorithm~\ref{alg:NI} 
        \State $\Mtil{u}^{(k+1)} \gets \Mtil{u}^{(k)} + 
                            \vek{\Delta}_{Z}(\Mtil{u}^{(k)})$
        \State $k\leftarrow k+1$
    \EndWhile 
  \end{algorithmic}
\end{algorithm}

The sequence generated by Algorithm~\ref{alg:BJNI} has been labelled
with a tilde $\{\Mtil{u}^{(k)}\}_k$ to distinguish it from the
exact sequence $\{\mat{u}^{(k)}\}_k$ generated by
Algorithm~\ref{alg:bJ}.   The question arises as to how well the
original sequence $\{\mat{u}^{(k)}\}_k$ is approximated by the one
produced non-intrusively by numerical integration
$\{\Mtil{u}^{(k)}\}_k$, and what its convergence behaviour is.
To that effect we partially cite and conclude from Theorem~4.1 in
\cite{Zander10}:
\begin{thm} \label{thm:thm4_1-Zander10} 
  Assume that the numerical integration in Algorithm~\ref{alg:NI} is
  performed such that 
  $\| \vek{G}_{\C{Q}}(\psi_{\alpha}P^{-1}R(\cdot;\tilde{u}^{(k)}))-
  \Delta_{Z, \alpha} \tilde{u}^{(k)}\|_{\C{U}} \leq \vepsilon / \sqrt{M} $,
  then the error in \refeq{eq:dqnumint} is estimated by 
  \begin{equation}  \label{eq:diff-NIBJ}
   \|\vek{\Delta}_{\C{Q}}(\Mtil{u}^{(k)}) -
   \vek{\Delta}_{Z}(\Mtil{u}^{(k)})  \|_{\C{U}^{\C{I}}}  \leq \vepsilon,
  \end{equation}
  and we have the following a posteriori error estimate for the iterates
  \begin{equation}  \label{eq:apost-NIBJ}
   \|\mat{u}^* - \Mtil{u}^{(k+1)}  \|_{\C{U}^{\C{I}}}  \leq 
   \frac{\vrho^*}{1-\vrho^*}
   \|\vek{\Delta}_{Z}(\Mtil{u}^{(k)})  \|_{\C{U}^{\C{I}}} +
  \frac{\vepsilon}{1-\vrho^*}.
  \end{equation}
In addition, we have that
  \begin{equation}  \label{eq:limit-NIBJ}
   \limsup_{k\to\infty} \|\mat{u}^* - \Mtil{u}^{(k)}  \|_{\C{U}^{\C{I}}}  \leq 
  \frac{\vepsilon}{1-\vrho^*}.
  \end{equation}
The satisfaction of the parametric equation may be gauged by
$ \|R(\cdot;\tilde{u}^{(k)})\|_{\C{Q}\otimes\C{U}}$.
\end{thm}
\begin{proof}
The \refeq{eq:diff-NIBJ} is a simple consequence of the assumption by
squaring and summing $M=|\C{I}|$ terms of size less than $\vepsilon /\sqrt{M}$,
and then taking the square root.
Everything else are then statements of Theorem~4.1 in \cite{Zander10}.
\end{proof}
The \refeq{eq:apost-NIBJ} shows that the modified sequence
$\{\Mtil{u}^{(k)}\}_k$ will not necessarily converge to
$\mat{u}^*$, even if $\vek{\Delta}_{Z}(\Mtil{u}^{(k)}) \to 0$
as $k\to\infty$, but \refeq{eq:limit-NIBJ} shows that it clusters
around $\mat{u}^*$ in a small neighbourhood.

\subsection{Computational effort and possible improvements}
   \label{SS:CompEff}
To assess the effort involved in a computational procedure and hence
its efficiency is always difficult, not least because it is not
always clear on how to measure computational effort.  Here we take
the view that the effort is only counted in solver calls, i.e.\
invocations of $S(p;u,R(p;u))$ \refeq{eq:It} or equivalently of
$P^{-1}R(p;u)$ \refeq{eq:It2}.  This means that the additional linear
algebra and computation of $\psi_\alpha(p_z)$ involved in the
Algorithms~\ref{alg:NI} and~\ref{alg:BJNI} is considered insignificant
in comparison to an invocation of the solver $S$.

The main contender for the Galerkin procedure outlined so far is to be
seen in what is called in the introduction a pseudo-spectral or
discrete projection, or a regression.  This can be described very quickly.
With $\{\psi_\alpha \}_{\alpha\in\C{I}}$ orthonormal, the coefficients
in the projection $u_{\C{I}}= \sum_{\alpha\in\C{I}} u^\alpha \psi_\alpha$ 
can be simply computed by inner products: 
\begin{equation}  \label{eq:direc-proj}
    u^\alpha = \ip{\psi_\alpha}{u^*}_{\C{Q}} = \int_{\C{P}}
    \psi_\alpha(p) u^*(p)\, \mu(\di p) \approx \sum_{z=1}^Z w_z \,
    \psi(p_z) u^*(p_z).
\end{equation}
One may remind oneself that this --- being the orthogonal projection
onto the subspace $\C{Q}_{\C{I}} \subseteq \C{Q}$ --- has the smallest
error to $u^*(p)$ in the norm $\| \cdot \|_{\C{Q}}$, but it does not
at all take into account the parametric equation.  The Galerkin
projection on the other hand will produce an approximation which is
optimal in minimising the residuum.  The approximation in
Algorithm~\ref{alg:proj} to \refeq{eq:direc-proj} is computed very
similarly to Algorithm~\ref{alg:NI}.
\begin{algorithm}
  \caption{Discrete projection according to \refeq{eq:direc-proj}}
    \label{alg:proj}
   \begin{algorithmic}
     \For {$\alpha \in \C{I}$}  
          \State $u^{\alpha} \gets 0 $ 
     \EndFor
        \Statex\Comment{ \%comment: the loop over integration points\%} 
     \For {$z\gets 1,\ldots,Z$} 
        \State Compute $u(p_z)$ acording to Algorithm~\ref{alg:basic}.
        \State $r_z \gets w_z\, u(p_z)$  
        \For {$\alpha \in \C{I}$}
          \State $u^{\alpha} \gets u^{\alpha} + \psi_\alpha(p_z) \, r_z$ 
        \EndFor
     \EndFor
  \end{algorithmic}
\end{algorithm}

The iterations from \refeq{eq:It} in Algorithm~\ref{alg:basic} with
one solver call per iteration in Algorithm~\ref{alg:proj} are assumed
to be contractions with contraction factor at most $\vrho^*$.  Say
that an iteration with contraction factor of $\vrho^*$ needs $L$
iterations to converge to the desired accuracy.  The discrete
projection needs $L$ solver calls on $Z$ integration points each,
i.e.\ $L \times Z$ solver calls.

The block Jacobi variant of the coupled Galerkin system in
Algorithm~\ref{alg:BJNI} needs \emph{one} solver call on $Z$
integration points.  But as it converges also with contraction
factor $\vrho^*$ --- see Corollary~\ref{coro:Lip-s}, it also needs $L$
iterations, i.e.\ in total also $L \times Z$ solver calls.

We see that in this measure of effort --- solver calls --- the
discrete projection and the block Jacobi iteration of the Galerkin
system need the same effort for comparable accuracy; something that is
borne out also in the numerical example in \refS{iterate}.  In case the
iteration in \refeq{eq:It} is quadratically convergent, e.g.\ it is
Newton's method, then this can not be matched by the block Jacobi
method; it will usually only have linear convergence.  When looking at
the other computations apart from the count of solver calls, in both
algorithms integrals have to be approximated by quadrature formulas.
In the discrete projection this happens only once, whereas in the
block Jacobi this is done in every global iteration.

Block Jacobi is probably the simplest method for coupled systems,
however it can be considerably accelerated \cite{MMV99a,
  matthiesStf03}, this ranges from the simple \emph{Aitken}
acceleration over block \emph{Gauss-Seidel} to \emph{Quasi-Newton}
methods.  In case the iterations from \refeq{eq:It} converge only
linearly, these extensions can then produce an advantage for the
Galerkin solution and may need considerably fewer than $L$ iterations,
as ``convergence information'' is shared for different values of $p$
or $\alpha$, something which will not happen in the decoupled discrete
projection.  Even \emph{Newton's} method \cite{MMV99a} can be emulated
on the global Galerkin system, where the action of the inverse of the
derivative on a vector is approximated by finite differences and
non-intrusive solver calls.  This last procedure is even able to
maintain quadratic convergence in case the iterations from
\refeq{eq:It} are quadratically convergent themselves.  These issues
will be taken up and published elsewhere.

Another area where considerable saving of work is possible in the
Galerkin procedure are sparse or low-rank approximations.  They come
about when viewing the solution --- and also other parametric elements
--- as \emph{tensors}, which may be used computationally in
\emph{low-rank} representations / approximations, see for example
\cite{Zander10}.  Again this is beyond the scope of the present paper,
and will published elsewhere.  Using such low-rank representations in
the originally uncoupled discrete projection produces a coupled
system, which then differs not substantially from the Galerkin system.
Other ways of building a low-rank representation were
already discussed in the introduction, and will be the subject of a
future paper.

%
%
%
%
%
%
%
%
%
%


%% file: chapters/iteration.tex
%

\section{Numerical example} \label{S:iterate}
\begin{figure}[h!]
\centering
\includegraphics[width=0.8\textwidth]{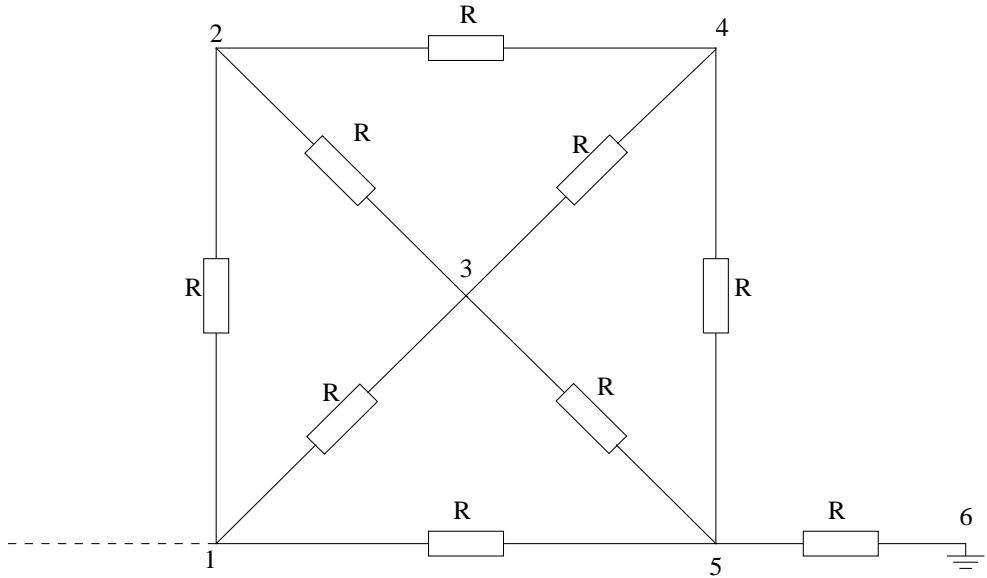}
\caption{Electrical resistor circuit.}
\label{F:circuit}
\end{figure}

Here we want to show the procedures discussed on a tiny example which
nonetheless is representative of parametric problems.  It is so simple
that it may be programmed with a few lines of code.
This computational example derives from a little electrical resistor
network with a global non-linearity.  The particular resistor network
we use is shown in \refig{circuit}.

Kirchhoff's and Ohm's laws lead to the following linear relation
between voltages $\vek{u}$ and currents $\vek{j}$ fed into the nodes,
where the numbering of the nodes corresponds to the equations --- 
node 6 is grounded ($u_6=0$) and so needs no equation, 
hence $\vek{u} \in \C{U} = \D{R}^5, \vek{K}\in \D{R}^{5\times 5}$:
\begin{equation} \label{eq:linSys}   
  \vek{K} \vek{u} = \vek{j}, 
\end{equation} 
with
\begin{equation} \label{eq:linMat}   
  \vek{K} = \frac{1}{R}\left [ \begin{array}{rrrrr} 
                                   3 & -1 & -1 & 0    &-1 \\
                                  -1 & 3 & -1 &-1    & 0   \\
                                  -1 & -1 & 4 &-1    & -1   \\
                                    0 & -1 & -1 &3    & -1   \\
                                  -1 &      0 & -1 &-1& 4      
      \end{array}   \right] , \;\;\;
\end{equation}
where we take all resistors to have the value $R=100$.

To make this toy system non-linear, we add a global cubic
non-linearity with an uncertain coefficient $(p_1+2) (\vek{u}^T
\vek{u})\, \vek{u}$.  We also make the feed-in current $\vek{f}= (p_2
+ 25)\vek{f}_0$ uncertain, so that we are eventually left to solve
this system for $\vek{u}$ (this is a concrete example of
\refeq{eq:I}):
\begin{align} \label{eq:lilSys}   
  \vek{A}(\vek{p};\vek{u}) :=& (\vek{K} \vek{u} +(p_1+2) (\vek{u}^T
  \vek{u})\, \vek{u}) = (p_2+25) \vek{f}_0 =: \vek{f}(\vek{p}),  \\
  \label{eq:lilSysf0}
       \vek{f}_0 :=& \left[1,  0, 0,  0, 0 \right]^T,
\end{align} 
where the random parameters $\vek{p}=(p_1,p_2)$ are assumed uniformly
and independently distributed in $[-1,1]$, and therefore we have for
the residuum (compare \refeq{eq:Res})
\begin{equation}  \label{eq:lilex_res}
  \vek{R}(\vek{p};\vek{u}) = (p_2+25) \vek{f}_0 - (\vek{K} \vek{u}
  +(p_1+2) (\vek{u}^T \vek{u})\, \vek{u}), 
\end{equation}
and for the preconditioner we take simply $\vek{P}=\vek{K} =
\Di_{\vek{u}} \vek{A}(\vek{p};\vek{0})$.

 The system can be solved in an iterative way as formulated in
 \refeq{eq:It}, with effectively
\begin{multline} \label{eq:iter_1}   
   \vek{u}^{(k+1)} = \vek{S}(\vek{p}; \vek{u}^{(k)},\vek{R}(\vek{p}; \vek{u}^{(k)}))
    = \vek{u}^{(k)} + \vek{P}^{-1}\vek{R}(\vek{p}; \vek{u}^{(k)}) = \\ 
   \vek{K}^{-1}  \left( (p_2+25) \vek{f}_0 - (p_1+2)
     ((\vek{u}^{(k)})^T  \vek{u}^{(k)}) \, \vek{u}^{(k)}
   \right)
 \end{multline} 
 The simple iteration \refeq{eq:iter_1} --- indeed a linearly
 convergent modified Newton method --- converges quite well for the
 chosen parameters.

 For the ansatz functions we take tensor products of Legendre
 polynomials, as they are orthogonal for the uniform measure, i.e.\
 we take $\psi_\alpha(\vek{p}) = \tilde{L}_\alpha(\vek{p}) =
 L_\alpha(\vek{p}) / \| L_\alpha \|$, the multi-variate normalised
 Legendre polynomial, and $L_\alpha(\vek{p}) = \prod_{i=1}^2
 \ell_{\alpha_i}(p_i)$, where the $\ell_i$ are the normal univariate
 Legendre polynomials, $\| L_\alpha \|= 4(2\alpha_1+1)^{-1}
 (2\alpha_2+1)^{-1}$, and $\alpha \in (\D{N}_0)^2$:
 \begin{align}  \label{eq:expansion}
 \vek{u}(\vek{p})  \approx \sum_{| \alpha|_1 \leq m}  \vek{u}_{\alpha}
 \tilde{L}_{\alpha}(\vek{p}) =: \vek{u}_{\C{I}}(\vek{p}),\quad \text{with}
 \end{align}
\[ \vek{u}_{\alpha} =  [ u_{\alpha,1},  \cdots,   u_{\alpha,5}]^T \in
\C{U} = \D{R}^5, \quad \text{and}\]
\[  \C{I} =\{\alpha = (\alpha_1, \alpha_2)\, :
\, |\alpha|_1 = \alpha_1 + \alpha_2 \le m\}   \subset (\D{N}_0)^2,
\;\;\; m \in \D{N};
 \]
hence for different $m\in \D{N}$ we will have different approximation
orders by polynomials of total degree $m$.

For the purpose of comparison we use two approaches to determine the
coefficients $\vek{u}_{\alpha}$ in \refeq{eq:expansion}, these are the
Galerkin approach according to Algorithm~\ref{alg:BJNI} with
numerically integrated residuum according to Algorithm~\ref{alg:NI}
for $\vek{u}_G(\vek{p})$, and collocation or more specifically discrete
projection with numerical integration according to
Algorithm~\ref{alg:proj} for $\vek{u}_C(\vek{p})$, both with the same
integration rule.  We choose here --- as we are only in two dimensions
--- a tensor-product Gauss-Legendre quadrature.  The quadrature order
was always taken so that products of test- and ansatz-functions
$\psi_\alpha \psi_\beta$ were integrated exactly for the chosen total
polynomial degree $m$ in \refeq{eq:expansion}.

First we computed a $N=1000$ sample Monte Carlo solution on random
points $\vek{p}_n \in \C{P}=[-1,1]^2, n=1,\ldots,N$ to high accuracy by
setting the convergence criterion in Algorithm~\ref{alg:basic} to
the machine tolerance.  These results were taken as the reference
solution for the following error estimation.  We computed the
root-mean-squared-error (RMSE) --- effectively the $L_2$ norm in
$ \C{Q} \otimes\C{U} \cong L_2([-1,1]^2;\D{R}^5)$ --- as 
\begin{equation}  \label{eq:RMSE}
  \epsilon_F = \left(\frac{1}{N} \sum_{n=1}^N \| \vek{u}(\vek{p}_n) -
    \vek{u}_F(\vek{p}_n)\|^2\right)^{1/2},
\end{equation}
where the \emph{functional approximation} method $F$ is either $G$
for the \emph{Galerkin} method or $C$ for the \emph{collocation}
method.
 
\begin{table}[htbp]
\begin{center}
\setlength{\extrarowheight}{7pt} 
 \begin{tabular}{c|c|c|c|c|c }
  \hline   \hline 
      \multirow{2}{2.5cm}{ Order of polynomial }  & $\epsilon_{tol}$ &
      \multicolumn{2}{c|}{\# of  solver evaluations }  &
      \multicolumn{2}{c}{RMSE $\epsilon_F$ from \refeq{eq:RMSE}} \\ 
  \cline{3-4}
     &  &Collocation &Galerkin & $F=C$ Collocation& $F=G$ Galerkin   \\
 \hline \hline
    m=2         & $10^{-6}$ & 73 &  81&$8.5\times 10^{-6}$ &
    $8.2\times 10^{-6}$  \\  
   \hline
   m=3          &$10^{-7}$  & 151 &  160 & $6.4\times 10^{-7}$ &
   $6.0\times 10^{-7}$\\    
      \hline
   m=4          &$10^{-8}$  & 268 & 300 & $4.2\times 10^{-8}$ &
   $4.0\times 10^{-8}$\\    
      \hline
   m=5          & $10^{-9}$& 430 & 468&$3.1\times 10^{-9}$ &
   $3.0\times 10^{-9}$\\        
  \hline   \hline  
\end{tabular}
 \end{center}
\caption{RMSE and number of solver evaluations  of collocation and
  Galerkin approaches}
\label{table:RMSE} 
\end{table}

The two approaches were carried out to compute the coefficients
$\vek{u}_\alpha$ in \refeq{eq:expansion}.  The criteria of convergence
for the iterative solvers were that the increment of $\vek u$ or $\vek
u_{\vek \alpha} $ is smaller than $\epsilon_{tol}$.  Tabulated in
\reftab{table:RMSE} are the $\epsilon_{tol}$ values obtained in a
sensitivity-range investigation such that further reduction of these
values would not improve accuracy, depending on the total polynomial
degree $m$.  The errors $\epsilon_F$ of each approach were estimated
as in \refeq{eq:RMSE} and are displayed in \reftab{table:RMSE},
together with the number of solver ($\vek{S}(\vek{p}, \vek{u} )$)
evaluations for total polynomial degrees $m=2, 3, 4$ and $5$.
The coefficients computed by either Galerkin or collocation differed
only in the eighth or ninth digit. 

It is seen in the results that in terms of ``the best possible
accuracy'' the Galerkin approach is slightly better than the
collocation one, though the former needs slightly more evaluations of
the solver.  This comes because the same convergence tolerance is used
for different equations in the two approaches.  The results essentially
confirm the theoretical analysis in \refS{galerkin}; for the same
accuracy both approaches need about the same number of solver calls,
i.e.\ the simple block Jacobi iteration of the Galerkin system
converges at essentially the same speed as the original iteration.

%
%
%
%
%
%
%
%
%


%% file: chapters/conclusion.tex
%

\section{Conclusion} \label{S:concl}
After reviewing the literature on numerical methods for parametric
equations, with a special emphasis on the subclass of stochastic
equations, we have introduced a general methodology to formulate
numerical methods relying on functional or spectral approximations.
We have shown that the Galerkin orthogonality conditions for the
residuum and the iteration equation are equivalent under certain
conditions, and that the simplest iterative scheme for the coupled
Galerkin system, the block Jacobi method, converges essentially at the
same speed as the original solver for a single parameter value.

In the main part for this ``plain vanilla'' Galerkin formulation, we
have shown how to approximate the preconditioned residuum in the
Galerkin equation through numerical integration, and the effects of
this on the iteration sequence.  Then these explicit non-intrusive
Galerkin algorithms have been compared on one simple, easy to
understand, baby example.  This showed that the theoretical analysis
was validated with these computations, and that even in the simplest
case of block Jacobi the Galerkin formulation is competitive with
collocation.  We finally recall once more the discussion in
\refS{resid} on possibilities to accelerate the coupled Galerkin
solution, something that is not possible for the decoupled collocation
approach.

%
%
%
%
%
%
